\documentclass[12pt]{article}
\usepackage[dvips]{graphicx}
\usepackage{amsmath, amssymb}
\usepackage{amsthm}

\newtheorem{theorem}{Theorem}[section]

\newtheorem{corollary}{Corollary}[section]
\newtheorem{remark}{Remark}[section]
\newtheorem{proposition}{Proposition}[section]

\newcommand{\pweight}{{Q}} % weights for paths
\newcommand{\relw}{\nu}    % relative betting weight
\newcommand{\rnp}{\rho}    % risk neutral probability

\newlength{\IndentI}
\newlength{\IndentII}
\newlength{\IndentIII}
\newlength{\WidthI}
\newlength{\WidthII}
\newlength{\WidthIII}
\title{Capital process and optimality properties of\\
 a Bayesian Skeptic in coin-tossing games} 

\author{
  Masayuki Kumon\\
  Risk Analysis Research Center\\
  Institute of Statistical Mathematics\\
  Akimichi Takemura\\
  Graduate School of Information Science and Technology\\
  University of Tokyo\\
  and\\
  Kei Takeuchi\\
  Faculty of International Studies\\
  Meiji Gakuin University }

\date{}

\begin{document}
\maketitle

\begin{abstract}
  We study capital process behavior in the fair-coin and
  biased-coin games in the framework of the game-theoretic probability
  of Shafer and Vovk (2001).  We show that if Skeptic uses a Bayesian
  strategy with a beta prior, the capital process is lucidly 
  expressed in terms of the past average of Reality's moves.  From
  this it is proved that the Skeptic's Bayesian strategy weakly forces
  the strong law of large numbers (SLLN) with the convergence rate of
  $O(\sqrt{\log n/n})$ and if Reality violates SLLN then the
  exponential growth rate of the capital process is very accurately
  described in terms of the Kullback divergence between the average of
  Reality's moves when she violates SLLN and the average when she
  observes SLLN.  We also investigate optimality properties
  associated with Bayesian strategy.
%   including the choice of the prior
%   for the Skeptic.
\end{abstract}

\noindent
{\it Keywords and phrases:} \ 
Azuma-Hoeffding-Bennett inequality,
beta-binomial distribution,
exchangeability,
game-theoretic probability, 
hypergeometric distribution,
Kullback divergence,
prior distribution.
%large deviation.

\section{Introduction}
\label{sec:intro}

Coin tossing is the most basic object in the study of probability theory.
This is also true in the new field of game-theoretic probability and
finance established by Shafer and Vovk (2001).  In fact they start with
discussing the fair-coin game in Section 3.1 of their book.  Recently
Kumon and Takemura (2007)\cite{kumon-takemura}, motivated by Takeuchi's works
(\cite{takeuchi:2004}, \cite{takeuchi:2004b}), 
proved that a very simple single strategy, based
only on the past average of Reality's moves, is weakly forcing SLLN with
the convergence rate of $O(\sqrt{\log n/n})$, which is a substantial
improvement over the original strategy of Shafer and Vovk.  
Versions of SLLN for unbounded moves by Reality is obtained in
Kumon, Takemura and Takeuchi (2007)\cite{ktt-unbounded}.

In this paper for general biased-coin games, we consider a class of Bayesian
strategies for Skeptic.  
As in
Kumon and Takemura (2007) we prove that Bayesian strategies in the class
weakly force SLLN with the convergence rate of $O(\sqrt{\log n/n})$.
Furthermore we establish the important fact that if Skeptic uses a Bayesian
strategy and Reality violates SLLN, then the exponential growth rate of
the Skeptic's capital process is very accurately described in terms of the
Kullback divergence between the average of Reality's moves when she
violates SLLN and the average when she observes SLLN.

In the protocol of the coin-tossing game of Shafer and Vovk (2001), there
is no probabilistic assumption on the behavior of Reality.  In the games,
Skeptic tries to become rich and Reality tries to prevent it.  However in
the Bayesian strategy, Skeptic simply and naively assumes that Reality
behaves probabilistically and for choosing his moves Skeptic uses the
Bayesian prediction of Reality's moves.  It is a remarkable fact that this
naive Bayesian prediction by Skeptic actually works and forces SLLN even
if Reality's moves are not probabilistic at all and Reality tries to beat
Skeptic as an adversary.  Furthermore Skeptic achieves an optimal growth
rate if Reality violates SLLN in a way accounted for by the prior.  As in
the standard statistical decision theory (e.g.\ Berger \cite{berger} and
Robert \cite{robert}), this optimality is inherent in Bayesian procedures.
However in the setting of the present paper, only a very simple protocol
of the game is assumed and no other modeling assumptions are made on
Reality's moves.  In this sense, we believe that the optimality considered
in this paper has much broader conceptual implications than those offered
by the standard statistical decision theory.

The organization of this paper is as follows.  In Section
\ref{sec:preliminary} we formulate coin-tossing games and set up some
notations.  In Section \ref{sec:bayesian} we formulate Bayesian strategy of
Skeptic for a given probability distribution on the space of Reality's
moves.  In  particular we consider Bayes procedures with respect to
beta-binomial distribution and hypergeometric distribution.
In Section \ref{sec:process} we show that Skeptic's capital process 
is expressed in a closed form when he uses a Bayesian strategy of Section
\ref{sec:bayesian}.  Then using Stirling's formula we prove that
the exponential growth rate of the capital process is naturally 
described by means of the Kullback divergence. 
% There we also give a simplified proof of Theorems 4.1 and 4.2 in 
% Kumon and Takemura (2005) \cite{kumon-takemura}. 
% In Section \ref{sec:choise-of-priors} we discuss choice of priors
% with some numerical studies and 
In Section \ref{sec:one-sided} we
consider the case that Skeptic is restricted to only buy the tickets in
coin-tossing games. 
In Section \ref{sec:numerical} some numerical examples are presented. 
We end the paper with some concluding remarks in Section \ref{sec:remarks}.

\section{Notations on coin-tossing game}
\label{sec:preliminary}

Here we prepare some notations on coin-tossing games.
We consider a general biased-coin game between Skeptic and Reality
in the following parametrization. In the protocol, the head and the tail of coin are coded as 1 and 0, respectively. Furthermore the probability of heads $0 < \rho < 1$ is given.

\medskip\noindent
\textsc{Biased-Coin Game} \\
\textbf{Protocol:}

\parshape=6
\IndentI   \WidthI
\IndentI   \WidthI
\IndentII  \WidthII
\IndentII  \WidthII
\IndentII  \WidthII
\IndentI   \WidthI
\noindent
${\cal K}_0 :=1$.\\
FOR  $n=1, 2, \dots$:\\
  Skeptic announces $M_n\in{\mathbb R}$.\\
  Reality announces $x_n\in \{-\rho,1-\rho\}$.\\
  ${\cal K}_n = {\cal K}_{n-1} + M_n x_n$.\\
END FOR

\medskip
\noindent

Skeptic starts with the initial capital ${\cal K}_0=1$. For each round
$n$, Skeptic first announces $M_n$, which is the number of tickets he
buys. $x_n$ stands for a ticket which sells for the price of $\rho$ and it pays the amount of one if Reality chooses heads and nothing if she chooses tails. $\rho$ is called the risk
neutral probability (e.g.\ Takemura and Suzuki (2007)).  
Except for Section \ref{sec:one-sided} we consider the case that Skeptic
is allowed to sell the tickets ($M_n<0$).  $M_n x_n$ is the increment of
Skeptic's capital from round $n-1$ to $n$.  The case $\rho = 1/2$ is
the fair-coin game.  Although most of our results can be presented in the
fair-coin game, it is instructive to consider the biased-coin game for the
additional asymmetry.

Let $\xi^n=x_1 \dots x_n$ denote the sequence of Reality's moves up to
time $n$.
$h_n$ denotes the number of heads and $t_n=n-h_n$ denotes the number
of tails in $\xi^n$. Furthermore write
\[
s_n = n \bar x_n = x_1 + \dots + x_n, \qquad s_0=0.
\]
Then $s_n = h_n - \rho n$ and 
\begin{equation}
\label{eq:heads-tails}
h_n=s_n + \rho n, \quad  \frac{h_n}{n}=\rnp + \bar x_n,
 \quad
t_n=n - (s_n + \rho n), \quad \frac{t_n}{n}= 1-\rnp - \bar x_n.
\end{equation}
In the above biased-coin game, as discussed below, Reality is forced to
observe SLLN $\bar x_n \rightarrow 0$ or equivalently to follow the risk
neutral probability $h_n/n \rightarrow \rnp$.

In the following we write 
\[\relw_n = \frac{M_n}{{\cal  K}_{n-1}}
\]
and express the recursive relation of the capital process in
the multiplicative form
$
{\cal K}_n = {\cal K}_{n-1} (1+ \relw_n x_n)
$. 
Then starting from the initial capital ${\cal K}_0=1$, the capital
process is written as
\[
{\cal K}_n = \prod_{i=1}^n (1+ \relw_i x_i).
%, \qquad \log {\cal K}_n = \sum_{i=1}^n \log (1+ \relw_i x_i).
\]
As in Shafer and Vovk (2001), we can think that Skeptic chooses his
strategy ${\cal P}$ which specifies
$\relw_i=\relw_i(x_1,\dots,x_{i-1})$ as a function of 
$x_1,\dots,x_{i-1}$.
Note that Skeptic has to choose 
\begin{equation}
\label{eq:restriction-on-nu}
 - \frac{1}{1-\rho} < \relw_i < \frac{1}{\rho}
\end{equation}
if he has to avoid bankruptcy (${\cal K}_i \le 0$).

In the protocol above there is no probabilistic assumption on the behavior
of Reality.  Suppose however that Skeptic models Reality's moves
$x_1,\dots, x_n$ by a probability distribution $\pweight=\pweight_n$ on
$\{-\rho,1-\rho\}^n$.  We denote the conditional expectation of $x_i$ under $\pweight$ 
given $x_1,\dots,x_{i-1}$ by
\[
\hat x_i^{\pweight} =  E_{\pweight}(x_i \mid x_1,\dots,x_{i-1}) = 
\hat p_i^\pweight - \rho,
\qquad i=1,\dots,n,
\]
where
\[
\hat p_i^\pweight 
= \hat p_i^\pweight (x_1,\dots,x_{i-1})=\pweight(x_i = 1-\rho \mid x_1,\dots,x_{i-1})
\]
is the conditional probability of heads.

\section{Bayesian strategy and its optimality}
\label{sec:bayesian}

In this section, we first discuss optimality inherent in a Bayesian
strategy.  It is not obvious how to formulate optimality of Skeptic's  strategies.  Allowing equalities in (\ref{eq:restriction-on-nu}), for each
fixed path of Reality's moves $\xi^n = x_1 \dots x_n$ up to time $n$, the
optimum moves of Skeptic are given as
\begin{equation}
\label{eq:meaningless}
\relw_i  = 
\begin{cases}
1/\rho & \text{if}\quad  x_i=1-\rho \\
 -1/(1-\rho) & \text{if}\quad  x_i=-\rho
\end{cases}
\end{equation}
for $i=1,\dots,n$, with the resulting capital 
${\cal K}_n=\rho^{-h_n}(1-\rho)^{-t_n}$.
However this argument is clearly meaningless because in each round 
Skeptic has to choose $\relw_i$ first and Reality chooses $x_i$ after
seeing $\relw_i$. For Skeptic's strategy ${\cal P}$ and a path $\xi^n$, we denote his log capital by $\log {\cal K}_n^{\cal P}(\xi^n)$ and denote a weight of the path by $w(\xi^n) \ge 0$. Then we consider the  weighted average of $\log {\cal K}_n^{\cal P}(\xi^n)$ and evaluate Skeptic's strategy ${\cal P}$ by
\begin{equation}
\label{eq:weighted0}
\sum_{\xi^n \in \{-\rho,1-\rho\}^n} w(\xi^n) \log {\cal K}_n^{\cal P}(\xi^n).
\end{equation}
There is no loss of generality in assuming 
$1=\sum_{\xi^n \in \{-\rho,1-\rho\}^n} w(\xi^n)$.  Therefore $\{w(\xi^n)\}$ is a
probability measure $\pweight$ and (\ref{eq:weighted0}) can be written as
\begin{equation}
\label{eq:weighted1}
E_\pweight[\log {\cal K}_n^{\cal P}(\xi^n)].
\end{equation}
For a given $\pweight$ the optimum ${\cal P}={\cal P}_Q$ maximizing
(\ref{eq:weighted1}) is given as follows by a standard backward
induction argument of dynamic programming.  We call ${\cal P}_Q$ the
Bayesian strategy for $Q$.

\begin{theorem}
\label{thm:bayes-optimum}
%Given a probability distribution  $\pweight$ on $\{-1,r-1\}^n$
% such that $0 < \hat p_i^\pweight (x_1,\dots, \allowbreak x_{i-1}) < 1$, 
% for all $i=1,\dots,n$ and  for all $x_1,\dots,x_n$.   Then 
The optimum strategy ${\cal P}={\cal P}_Q$ maximizing 
(\ref{eq:weighted1}) is given by
\begin{equation}
\label{eq:optimum-move}
\relw_i(x_1,\dots,x_{i-1})=\frac{\hat p_i^Q - \rho}{\rho(1-\rho)}
= \frac{\hat x_i^Q}{\rho(1-\rho)}
=\frac{E_{\pweight}(x_i \mid x_1,\dots,x_{i-1})}{\rho(1-\rho)}, \quad i=1,\dots,n.
\end{equation}
\end{theorem}

\begin{proof}
Consider the optimum move of Skeptic at the last round $n$, given all
previous moves of the players.  The terms of 
\[
\log {\cal K}_n^{\cal P}(\xi^n) = \sum_{i=1}^n
\log (1+\relw_i(x_1,\dots,x_{i-1})x_i)
\]
are all fixed except for the last term $\log(1+\relw_n x_n)$.  Consider
maximizing the conditional expected value 
$g(\nu)=E_\pweight( \log (1+\nu x_n) \mid x_1,\dots,x_{n-1})$.
$g(\nu)$ is a concave function of $\nu$ and under the convention $0=0
\log 0$, $g(\nu)$ is maximized by
solving $g'(\nu)=0$  as long as this equation has a solution in 
$[- 1/(1-\rho),1/\rho]$.  Now
\[
g'(\nu)= \hat p_n^\pweight \frac{1-\rho}{1+\nu (1-\rho)} -
(1-\hat p_n^\pweight) \frac{\rho}{1-\nu \rho}.
\]
Solving this we have
\[
\relw_n(x_1,\dots,x_{n-1})= \frac{\hat p_n^\pweight-\rho}{\rho(1-\rho)} = \frac{\hat x_n^\pweight}{\rho(1-\rho)}.
\]
Note that  $(\hat p_n^\pweight-\rho)/\rho(1-\rho) \in [- 1/(1-\rho),1/\rho]$.  Therefore the optimum
move of the Skeptic at the last round $i=n$ is given by 
(\ref{eq:optimum-move}).  Now for the round $i=n-1$ we can do the same
argument based on the conditional distribution of $x_{n-1}$ given
$x_1,\dots,x_{n-2}$.  Then by backward induction (\ref{eq:optimum-move})
is proved for all $i=n-1,\dots,1$.
\end{proof}

In the beginning of this section we argued that (\ref{eq:meaningless}) for
a fixed sequence $\xi^n=x_1\dots x_{n}$ is meaningless.  However Theorem
\ref{thm:bayes-optimum} reduces to (\ref{eq:meaningless}) if $Q$ is a
point mass at a particular $\xi^n$.  This shows that the optimality in
Theorem \ref{thm:bayes-optimum} for a given $Q$ is in a sense a trivial
inherent optimality of a Bayes procedure and the important question is the
choice of the probability measure $Q$.  In the following we consider
various choices of $Q$, including those whose support is a proper subset
of $\{ -\rho,1-\rho\}^n$.

The obvious first candidate for $Q$ is the Bernoulli model, where Reality
is modeled to choose her move independently for each round as
$
Q(x_i=1-\rho)=p=1-Q(x_i=-\rho).
$
Then
\[
Q(\xi^n)=p^{h_n}(1-p)^{t_n}, 
\]
where $h_n$ and $t_n$ are the numbers of heads and tails in $\xi^n = x_1 \dots
x_n$.  In this case $\hat p_i^\pweight=p$ for all $i$ and $\relw_i =
(p-\rho)/\rho(1-\rho)$ is a constant.  We see that the fixed $\epsilon$-strategy of
Chapter 3 in Shafer and Vovk (2001) corresponds to this case.

{}From Bayesian viewpoint it is natural to consider a mixture of the
Bernoulli model with respect to a beta prior, which is the natural conjugate
to the binomial distribution.  
Suppose that Skeptic has a prior
beta distribution $\pi$ on the success probability $p$ of the Bernoulli
model.
\begin{equation}
\label{eq:beta-prior}
\pi(dp) =  \frac{1}{B(\alpha,\beta)}  
p^{\alpha-1} (1-p)^{\beta -1}dp ,  \quad \alpha,\beta > 0.
\end{equation}
We can think of $\alpha$ and $\beta$ as prior numbers of heads and tails.
Then
\begin{align}
Q(\xi^n) 
&= \frac{1}{B(\alpha,\beta)}  \int_0^1 p^{h_n+\alpha-1} (1-p)^{t_n + \beta -1}dp
= \frac{B(\alpha+h_n, \beta+t_n)}{B(\alpha,\beta)}
\nonumber \\
&= \frac{(\Gamma(\alpha+h_n)/\Gamma(\alpha)) \times
         (\Gamma(\beta+t_n)/\Gamma(\beta))}
        {\Gamma(\alpha+\beta+n)/\Gamma(\alpha+\beta)}
\label{eq:beta-binomial-prob}
\end{align}
is the beta-binomial distribution 
and we call this model ``beta-binomial model'' and call the Skeptic's
associated strategy ``beta-binomial strategy''. In this model
\begin{equation}
\label{eq:beta-binomial}
\hat p_n^\pweight = \frac{B(\alpha+ h_{n-1}+1 , \beta+t_{n-1})}{B(\alpha+
  h_{n-1} , \beta+t_{n-1})} = \frac{\alpha+h_{n-1}}{\alpha+\beta+ n-1},
\end{equation}
and it follows that
\begin{equation}
\label{eq:rho-beta-binomial}
\relw_i =\frac{\hat p_i^\pweight-\rho}{\rho(1-\rho)}
%= \frac{r\alpha + s_{i-1} +i-1-(\alpha+\beta+i-1) }{(\alpha+\beta+i-1)(r-1)}
= \frac{(1-\rho)\alpha - \rho \beta + s_{i-1}}{(\alpha+\beta+i-1)\rho(1-\rho)} , \quad i=1,\dots,n.
\end{equation}

Finally we consider a hypergeometric model, whose support may be a proper
subset of $\{ -\rho,1-\rho\}^n$.  This is a somewhat artificial model, but it
is useful as a benchmark in considering optimality of the beta-binomial strategy.
% choice of priors in Section
%\ref{sec:choise-of-priors} as a benchmark.  
As we see below the
hypergeometric model has an opposite characteristic to the beta-binomial
model.  Let $N \ge n$ and $0 \le M \le N$ be positive
integers.  

Consider an urn containing $M$ red balls and $N-M$ black balls. Skeptic
models Reality's behavior in such a way that Reality draws a ball from the
urn one by one without replacement and she chooses heads if the ball is
red in the $i$-th drawing.  Therefore Skeptic is considering an initial
part of a finite horizon game with $N$ rounds.  Under this model
\begin{align}
\label{eq:hypergeometric-prob}
Q(\xi^n)
&=\frac{1}{\binom{n}{h_n}} \frac{\binom{M}{h_n}\binom{N-M}{t_n}}{\binom{N}{n}}
= \frac{(M!/(M-h_n)!)\times ((N-M)!/(N-M-t_n)!)}{N!/(N-n)!} , \\
& \hspace{5cm} 0 \le h_n \le M, \ 0 \le t_n \le N-M,
\nonumber
\end{align}
and 
\begin{equation}
\hat p_n^\pweight = \frac{Q(x_1 \ldots x_{n-1} (1-\rho))}{Q(x_1 \ldots x_{n-1})}
= \frac{(M-h_{n-1})!}{(M-h_{n-1}-1)!}
\frac{(N-n)!}{(N-n+1)!}
=\frac{M-h_{n-1}}{N-n+1}.
\label{eq:alpha-hypergeometric}
\end{equation}
Since the hypergeometric model holds for each $i\le n$, it follows that
\begin{equation}
\label{eq:rho-hypergeometric}
\relw_i 
= \frac{M - s_{i-1}-\rho(i-1) - \rho(N-i+1)}{(N-i+1)\rho(1-\rho)}
= \frac{M-\rho N - s_{i-1}}{(N-i+1)\rho(1-\rho)} .
\end{equation}
Note the similarities between
(\ref{eq:beta-binomial-prob})
and
(\ref{eq:hypergeometric-prob})
and between
(\ref{eq:beta-binomial}) and 
(\ref{eq:alpha-hypergeometric}).  
If we put
\begin{equation}
\label{eq:correspondence}
\alpha=-M, \quad \beta=-(N-M),
\end{equation}
then (\ref{eq:beta-binomial}) coincides with
(\ref{eq:alpha-hypergeometric}).  

In order to make the correspondence clearer, define
\[
(a)_m = a(a+1)\cdots (a+m-1), \qquad (a)_0=1, 
\]
for real $a$ and non-negative integer $m$.   For $a>0$, 
$(a)_m =\Gamma(a+m)/\Gamma(a)$.  Then with the correspondence
(\ref{eq:correspondence}),
both (\ref{eq:beta-binomial-prob}) and (\ref{eq:hypergeometric-prob}) can
be written as
\begin{equation}
\label{eq:qxin}
Q(\xi^n) = \frac{(\alpha)_{h_n} (\beta)_{t_n}}{(\alpha+\beta)_{n}}.
\end{equation}
We see that the hypergeometric model is the ``negative'' of the beta-binomial model.
% we obtain
% \[
% \frac{-Mr + N + s_{i-1}}{(-N+i-1)(r-1)}
% = \frac{rM - N -  s_{i-1}}{(N-i+1)(r-1)}
% \]
In the beta-binomial model $\hat p_n^\pweight$ is an increasing
function of $h_{n-1}$, whereas in the hypergeometric model $\hat
p_n^\pweight$ is a decreasing function.  Another way of understanding
the connection is the Polya's urn model (e.g.\ Section V.2 of Feller
(1968), Takemura and Suzuki (2005)). % \cite{takemura-suzuki}).
Hypergeometric model corresponds to sampling without replacement,
Bernoulli model corresponds to sampling with replacement and
beta-binomial model corresponds to Polya's urn model where a ball with
the same color is added at each draw.

The extreme case of the hypergeometric model is the case $N=n$.  Then $Q$
is concentrated on $\xi^n$ with the number of heads exactly equal to $M$.
Among the exchangeable models, in the sense that they assign the same
probability to all $\xi^n$ with the same number of heads, this model is
most ``meaningless'' as in the very first example of this section.  Note
that the support of the hypergeometric model is the whole sample space $\{-\rho,1-\rho\}^n$ if and only if $n \le \min(M,N-M)$. One way of justifying the hypergeometric model is to add a requirement in the protocol of the game,
such that Reality has to choose her path in the support of $Q$.  By this
requirement Reality's move space is decreased and the game becomes more
favorable to Skeptic.  This implies that the hypergeometric case, in
particular the extreme case $N=n$, can
serve as an upper bound to Skeptic's capital process.  We should also
mention de Finetti's theorem (\cite{diaconis-freedman},
\cite{takeuchi-takemura:1987}), which states that an infinite sequence of
exchangeable 0-1 random variables has to be a mixture of infinite
independent Bernoulli trials.  We see that the hypergeometric model is naturally
associated with a finite horizon game.

\section{Capital precess of Bayesian strategy} % and weak forcing of SLLN}
\label{sec:process}

In this section we investigate capital process when Skeptic uses a
Bayesian strategy.  We first give a general formula for Skeptic's
capital as a ratio of the probabilities of Reality's path under the
assumed distribution and under the risk neutral probability distribution.
This gives us a closed form expression of the capital process for the
beta-binomial model and the hypergeometric model of the previous
section.  Then for the beta-binomial model, we use Stirling's formula
to describe the exponential growth rate of the capital process. 
It leads to the Kullback divergence and proves that the beta-binomial
strategy weakly forces SLLN with the convergence rate of $O(\sqrt{\log
  n/n})$.  Finally we consider optimality of the beta-binomial strategy in relation to the hypergeometric model.

% We now show that for the Bayesian strategy of the previous section the
% value of the capital process at time $n$ only depends on $s_n$ and
% hence independent of the order of heads and tails up to time $n$.
% In the
% following write $c=\alpha+\beta$ and $d=(r-1)\alpha - \beta$ in
% By the correspondence (\ref{eq:correspondence}), the beta-binomial model
% and the hypergeometric model can be handled together.

The following theorem shows that the capital process for a Bayesian
strategy can be written as the ratio of probabilities under the
assumed model and under the risk neutral probability measure.

\begin{theorem}
\label{thm:path}
Let $Q=Q_n$ be a probability measure on the set of paths $\{\xi^n\}$
of length $n$ and let ${\cal P}_Q$ denote the Bayesian strategy for
$Q$.  The value of the capital process ${\cal K}_n^{{\cal P}_Q}$ for
${\cal P}_Q$ is given by
\begin{equation}
\label{eq:thm-2}
{\cal K}^{{\cal P}_Q}_n (\xi^n)= \frac{Q(\xi^n)}{\rnp^{h_n}
  (1-\rnp)^{t_n}} ,
% =
% \frac{1}{\rnp^{h_n} (1-\rnp)^{t_n}}
% \frac{(\alpha)_{h_n} (\beta)_{t_n}}{(\alpha+\beta)_{n}}
\end{equation}
where $0 < \rnp < 1$ is the risk neutral probability. % and $(1-\rnp)= 1-\rnp$.
\end{theorem}

\begin{proof}
We prove (\ref{eq:thm-2}) by induction on $n$.  Recall that
$\relw_i=(\hat p_i^Q -\rho)/\rho(1-\rho)$.  Consider $n=1$.  If $x_1=1-\rho$, then 
\[
{\cal K}_1^{{\cal P}_Q} = 1 + \relw_1 (1-\rho) = 1 + \frac{\hat p_1^Q -\rho}{\rho} = \frac{\hat p_1^Q}{\rho}=\frac{Q(x_1 = 1-\rho)}{\rho}.
\]
On the other hand if $x_1 = -\rho$, then
\[
{\cal K}_1^{{\cal P}_Q} = 1- \relw_1 \rho= \frac{1-\hat p_1^Q}{1-\rho}
= \frac{Q(x_1 = -\rho)}{1-\rho}.
\]
This proves (\ref{eq:thm-2}) for $n=1$.

Now suppose that (\ref{eq:thm-2}) holds up to $n-1$. Then
\[
{\cal K}^{{\cal P}_Q}_n (\xi^n)=\frac{Q(\xi^{n-1})}{\rho^{h_{n-1}}
  (1-\rho)^{t_{n-1}}} (1+\nu_n x_n).
\]
As in the case of $n=1$  it holds that
\[
1+\nu_n x_n= \begin{cases}
Q(x_n = 1-\rho \mid \xi^{n-1})/\rho & \text{if}\quad   x_n = 1-\rho\\
Q(x_n = -\rho \mid \xi^{n-1})/(1-\rho)  & \text{if}\quad   x_n = -\rho.
\end{cases}
\]
Therefore (\ref{eq:thm-2}) holds also for  $n$.
\end{proof}

\begin{corollary}  For the beta-binomial model and the hypergeometric
  model
\begin{equation}
\label{eq:thm-2a}
{\cal K}^{{\cal P}_Q}_n (\xi^n)= 
\frac{1}{\rnp^{h_n}
  (1-\rnp)^{t_n}} \frac{(\alpha)_{h_n} (\beta)_{t_n}}{(\alpha+\beta)_n}.
\end{equation}
\end{corollary}

\begin{remark}
\label{rem:vovk}
Formulation of Theorem \ref{thm:path} in the present form was
suggested by Vladimir Vovk to one of the authors in a discussion
during the 16th international conference on algorithmic learning
theory.  Theorem \ref{thm:path} actually follows from some general
facts, including the equivalence of game-theoretic martingales and
measure-theoretic martingales (Section 8.2 of Shafer and Vovk (2001)),
expressing positive martingales with expected value of 1 as likelihood
ratios and the non-negativeness of the Kullback divergence.
\end{remark}

We now use Stirling's formula to prove that the beta-binomial
strategy weakly forces SLLN.  Let $\alpha>0, \beta>0$ in
(\ref{eq:thm-2a}).  
% Stirling's formula
% \[
% \log \Gamma(x)=(x-1/2)\log x - x + \log\sqrt{2\pi}
% + \sum_{n=1}^\infty \frac{B_{2n}}{2n (2n-1) x^{2n-1}}
% \]
% $B_0=1, B_1 = -1/2, B_2 = 1/6, B_3 = 0, B_4 = -1/30$
The log capital process is written as
\begin{align*}
\log {\cal K}^{\cal P}_n (\xi^n)
&= -h_n \log \rnp - t_n \log (1-\rnp) 
+ \log \Gamma(\alpha+h_n) - \log \Gamma(\alpha)\\
&\qquad + \log \Gamma(\beta+t_n) - \log \Gamma(\beta)\\
&\qquad - \log\Gamma(\alpha+\beta+n) + \log \Gamma(\alpha+\beta).
\end{align*}
If both $h_n$ and $t_n$ are large, we can use Stirling's formula
\[
% \log \Gamma(x)=(x-\frac{1}{2})\log x - x + \log\sqrt{2\pi} + \frac{1}{12x}
% + O(x^{-3})
\log \Gamma(x)=\Bigl(x-\frac{1}{2}\Bigr)
\log x - x + \log\sqrt{2\pi} + O(x^{-1}).
% + \sum_{n=1}^\infty \frac{B_{2n}}{2n (2n-1) x^{2n-1}}.
% \]
\]
More precisely  for all $x>0$
\[ 0 < 
\log \Gamma(x) - \Big[\Bigl(x-\frac{1}{2}\Bigr)
\log x - x + \log\sqrt{2\pi}\Big]
< \frac{1}{12x}.
\]

For notational simplicity write $n'=\alpha+\beta+n$,
$h_n'=\alpha+h_n$,  $t_n'=\beta + t_n$. Then % $n'=h_n'+t_n'$  and
\begin{align*}
& \log \Gamma(h_n')  + \log \Gamma(t_n') -
\log\Gamma(n') \\
&\qquad \qquad =
h_n'\log h_n' + t_n'\log t_n' - n' \log n' -\frac{1}{2}\log\frac{ h_n' t_n'}{n'}
+ \log\sqrt{2\pi}+ O\biggl(\frac{1}{\min (h_n', t_n')}\biggr)\\
&\qquad \qquad =
 h_n' \log \frac{h_n'}{n'} + t_n' \log \frac{t_n'}{n'} 
-\frac{1}{2}\log\frac{ h_n' t_n'}{n'}
+ \log\sqrt{2\pi}+ O\biggl(\frac{1}{\min (h_n', t_n')}\biggr).
\end{align*}
For $0< p,\, q < 1$, let
\[
% \begin{equation}
% \label{eq:kullback}
   D(p\|q) = p\log\frac{p}{q} + (1 - p)\log\frac{1 - p}{1 - q}
%\end{equation}
\]
denote the Kullback divergence between $p$ and $q $. Then 
%writing $h_n = h_n' - \alpha$ and $t_n = t_n' - \beta$, 
$\log {\cal K}_n^{\cal P}$
is written as
\begin{equation}
\label{eq:log capital}
\log {\cal K}_n^{\cal P}
= n' D\Bigl(\frac{h_n'}{n'} \Big\| \rnp\Bigr)  
-\frac{1}{2}\log\frac{ h_n' t_n'}{n'} +
c_0(\alpha,\beta) + O\biggl(\frac{1}{\min (h_n', t_n')}\biggr), 
\end{equation}
where
\[
c_0(\alpha,\beta)= -\log B(\alpha,\beta) + \alpha \log \rnp +  \beta \log
(1-\rnp) +\log\sqrt{2\pi}.
\]

By the Taylor expansion
\begin{align*} 
D(\rnp + \delta \| \rnp) = 
(\rnp + \delta)\log \Bigl(1 + \frac{\delta}{\rnp}\Bigr) + 
(1-\rnp - \delta)\log \Bigl(1 - \frac{\delta}{1-\rnp}\Bigr) 
= \frac{\delta^2}{2\rnp^*(1-\rnp^*)}, 
\end{align*}
where $\rnp^*$ is some value between $\rnp$ and $\rnp + \delta$.
Recall that $h_n/n=\rnp + \bar x_n$.
Then with $\delta = \bar{x}_n$ %\ (|\delta| < 1)$
we have
\begin{align}
\label{eq:log capital2}
\log {\cal K}^{\cal P}_n &= 
% \frac{\rnp}{2(1-\rnp)} 
% n\bar{x}_n^2 - \frac{1}{2}\log n + O(n\bar{x}_n^3) 
% = \frac{1}{2}\Bigl(\frac{\rnp}{1-\rnp}n\bar{x}_n^2 - \log n \Bigr) + 
% O(n\bar{x}_n^3).
\frac{n\bar x_n^2}{2\rnp^*(1-\rnp^*)} - \frac{1}{2} \log n + O(1)\nonumber\\
&= \frac{1}{2}\biggl(\frac{n\bar{x}_n^2}{\rho(1-\rnp)} - \log n \biggr) + o(n\bar x_n^2)+ O(1).
\end{align}
Hence we obtain the following result.

\begin{theorem}
\label{thm:rate} When Skeptic follows the beta-binomial strategy ($\alpha>0$, $\beta>0$) given by (\ref{eq:rho-beta-binomial}),
\begin{equation*}
\limsup_n \ \  (\liminf_n)\ \  {\cal K}^{\cal P}_n = \infty
\end{equation*}
if and only if 
\begin{equation*}
\limsup_n\ \  (\liminf_n)\ \ 
\biggl(\frac{n\bar{x}_n^2}{\rho(1-\rnp)} - \log n \biggr) = \infty.
\end{equation*}
A sufficient condition for $\limsup_n{\cal K}^{\cal P}_n = \infty$ is
\begin{equation*}
\limsup_n \frac{\sqrt{n}|\bar{x}_n|}{\sqrt{\log n}} > 
\sqrt{\rho(1-\rnp)},
\end{equation*}
and a necessary condition for $\limsup_n{\cal K}^{\cal P}_n = \infty$ is
\begin{equation*}
\limsup_n \frac{\sqrt{n}|\bar{x}_n|}{\sqrt{\log n}} \ge 
\sqrt{\rho(1-\rnp)}.
\end{equation*}
\end{theorem}

Note that $\lim {\cal K}_n^{\cal P}=\infty$ if and only if $\liminf {\cal K}_n^{\cal P}=\infty$.
This theorem states that the Bayesian strategy (\ref{eq:rho-beta-binomial}) 
weakly forces that $\bar{x}_n$ converges to 0. 
The convergence rate is $O(\sqrt{\log n/n})$ and the the convergence factor is $\sqrt{\rho(1-\rnp)}$.
%or the biasedness of coin.

We also note that from the log expression (\ref{eq:log capital2}), 
the capital ${\cal K}_n^{\cal P}$ behaves as
\begin{align*}
{\cal K}_n^{\cal P} \simeq n^A,\qquad 
A = \frac{1}{2\rho(1 - \rho)}\biggl(\frac{n\bar{x}_n^2}{\log n} - \rho(1 - \rho) \biggr).
\end{align*}
Hence we know that,  if SLLN holds ($\bar{x}_n \to 0$) with
a convergence rate slower than $O(\sqrt{\log n/n})$, then Skeptic's capital grows faster than any polynomial order of $n$.

\begin{remark}
\label{rem:4.1}
As clarified by the above argument, the Bayesian strategy which is
plain in itself, greatly simplifies the proof and derivation of the
strong law of large numbers for coin-tossing games. Furthermore Skeptic
has a wide choice of prior distributions on the behavior of Reality,
although the question of optimal choice of prior distributions seems
to be a difficult problem.
% This subject will be 
% pursued in the next section and also in 
Capital processes for various values of $\alpha, \beta$ are illustrated
by numerical examples in Section \ref{sec:numerical}.
\end{remark}

\begin{remark}
\label{rem:4.2} As discussed in Remark 2 of Kumon and Takemura
\cite{kumon-takemura},
the convergence rate of $O(\sqrt{\log n/n})$ in Theorem
\ref{thm:rate} is closely
related to Azuma-Hoeffding-Bennett inequality, which is a large
deviation type inequality for bounded martingale differences.  See
Azuma (1967), %, \cite{azuma},
Dembo and Zeitouni (1998), % \cite{dembo-zeitouni},
and Appendix A.7 of Vovk, Gammerman and Shafer
(2005). % \cite{vovk-gammerman-shafer-2004}.
\end{remark}

\begin{remark}
We also point the importance of game-theoretic results that the
question of the optimal growth rate of Skeptic's capital cannot even
be asked in the standard measure-theoretic probability theory. This
question can be asked in the algorithmic theory of randomness, see
e.g. \ Schnorr (1970) (1971), Vovk (1987). 
But game-theoretic probability theory does not suffer from the two serious disadvantages of the algorithmic theory of randomness: the arbitrary constants and the heavy restrictions on the allowed sample spaces caused by considerations of computability.
\end{remark}

For the rest of this section we consider performance of the beta-binomial
strategy compared to the hypergeometric model. As remarked at the end
of Section \ref{sec:bayesian}, the extreme case $N=n$ of the hypergeometric
model can serve as an upper bound for the beta-binomial models.  
We first state this in the following proposition.

\begin{proposition}  For every exchangeable probability $Q$ on the set
  of paths $\{\xi^n\}$  of length $n$,
  the following inequality holds for each $\xi^n$.
\begin{equation}
\label{eq:exchangeable-optimality}
%\frac{1}{\rnp^{h_n}(1-\rnp)^{t_n}} 
\frac{1}{\rho^{h_n}(1-\rho)^{t_n}}\frac{h_n! t_n!}{n!} \ge 
%\frac{1}{\rnp^{h_n} (1-\rnp)^{t_n}}
%\frac{(\alpha)_{h_n} (\beta)_{t_n}}{(\alpha+\beta)_{n}}
%\frac{B(\alpha+h_n, \beta+t_n)}{B(\alpha,\beta)}.
{\cal K}_n^{{\cal P}_Q}(\xi^n).
\end{equation}
%\ref{thm:bayes-optimum}
\end{proposition}

\begin{proof} 
The left-hand side is the value of the capital process for
the hypergeometric prior with $n=N$. The prior is the uniform
distribution over the set of paths $\{\xi^n\}$ with the same number
$h_n$ of heads.  For exchangeable $Q$, the right-hand side is
constant for each path $\xi^n$ with the same number $h_n$ of heads. Then the
inequality follows from the optimality of the left-hand side for the
hypergeometric prior with $n=N$.
\end{proof}

Write
\[
{\cal K}_n^*(\xi^n) = \frac{1}{\rnp^{h_n}(1-\rnp)^{t_n}}\frac{h_n! t_n!}{n!}. 
\]
% which is the value of the capital process for the extreme
% hypergeometric prior with $N=n$. 
Conceptually we need to consider
${\cal K}_n^*$ for each $n$ and the number of heads $h_n$, separately. 
Let ${\cal K}_n^{\alpha,\beta}$
denote the capital process for the beta-binomial model.  
Then ${\cal K}_n^*(\xi^n) /{\cal K}_n^{\alpha,\beta}(\xi)$ is 
written as
\[
\frac{{\cal K}_n^*(\xi^n)}{{\cal K}_n^{\alpha,\beta}(\xi^n)}
= \frac{h_n! t_n!/n!}{B(\alpha+h_n, \beta+t_n)} B(\alpha,\beta) \ge 1.
\]
Therefore
\begin{align*}
\log {\cal K}_n^*(\xi^n) - \log {\cal K}_n^{\alpha,\beta}(\xi^n)
&=\log h_n! + \log t_n! - \log n!  \\
&\quad  -\big(\log \Gamma(\alpha+h_n) + \log
\Gamma(\beta+t_n) - \log \Gamma(\alpha + \beta +n)\big)\\
&\quad   + \log B(\alpha,\beta).
\end{align*}
Stirling's formula for $\log x!$ is written as
\[
0 < \log x!  - 
\Big[\Bigl(x+\frac{1}{2}\Bigr)\log (x + 1) - (x + 1) + 
\log \sqrt{2\pi}\Big] < \frac{1}{12(x + 1)}.
% = \Bigl(x+\frac{1}{2}\Bigr)\log (x + 1) - (x + 1)  + 
% \log \sqrt{2\pi} +O((x + 1)^{-1}).
\]
Hence we have 
%$\log {\cal K}_n^*(\xi^n) - \log {\cal  K}_n^{\alpha,\beta}(\xi^n)$ is written as
\begin{align*}
\log {\cal K}_n^*(\xi^n) - \log {\cal K}_n^{\alpha,\beta}(\xi^n)
&=
\Bigl(h_n + \frac{1}{2}\Bigr) \log h_n 
 - \Bigl(h_n + \alpha - \frac{1}{2}\Bigr) \log (h_n+ \alpha) \\
&\quad + \Bigl(t_n + \frac{1}{2}\Bigr) \log t_n 
 - \Bigl(t_n + \beta - \frac{1}{2}\Bigr) \log (t_n + \beta) \\
&\quad - \Bigl(n + \frac{1}{2}\Bigr) \log n 
 + \Bigl(n + \alpha + \beta - \frac{1}{2}\Bigr) 
\log (n + \alpha + \beta) \\
&\quad + \log B(\alpha,\beta) + O\biggl(\frac{1}{\min (h_n',
  t_n')}\biggr) .
\end{align*}
The right-hand side is further simplified as
\begin{align*}
&(1-\alpha)\log h_n - \alpha + (1-\beta) \log t_n - \beta 
-(1-\alpha -\beta) \log n + (\alpha+\beta)\\
&\qquad + \log B(\alpha,\beta) + O\biggl(\frac{1}{\min (h_n', t_n')}\biggr)\\
&\ = \log n + (1-\alpha)\log\frac{h_n}{n} + (1-\beta)\log\frac{t_n}{n}
+  \log B(\alpha,\beta) + O\biggl(\frac{1}{\min (h_n', t_n')}\biggr).
\end{align*}
We summarize the above calculation in the following proposition, which states that ${\cal K}_n^*(\xi^n)$ surpasses 
${\cal K}_n^{\alpha,\beta}(\xi^n)$ only by a polynomial factor of $n$.
\begin{proposition}
\[
\frac{{\cal K}_n^*(\xi^n)}{{\cal K}_n^{\alpha,\beta}(\xi^n)}
= n B(\alpha,\beta) \Big(\frac{h_n}{n}\Big)^{1-\alpha}
  \Big(\frac{t_n}{n}\Big)^{1-\beta}
 (1+O(1/\min (h_n', t_n')).
\]
\end{proposition}
Considering the fact that Skeptic can achieve ${\cal K}_n^*(\xi^n)$ 
only in ``hindsight'' (i.e.\ after seeing the number of heads $h_n$ for
each $n$)  and the fact that
${\cal K}_n^*(\xi^n)$ and ${\cal K}_n^{\alpha,\beta}(\xi^n)$ grow
exponentially when Reality violates SLLN, 
we see that the beta-binomial strategy is close to optimum.

\section{Capital process for the one-sided case}
\label{sec:one-sided}
In this section we assume the protocol as before, but now Skeptic 
is required to choose $M_n \ge 0$. We consider Skeptic's strategy ${\cal P}^+$ with
\begin{align}
\label{eq:one-sided strategy}
\relw^+_i = \max(\relw_i, 0),
\end{align}
where $\relw_i$ is given by (\ref{eq:rho-beta-binomial}).  We can also
consider the negative part strategy ${\cal P}^-$ with 
$\relw^-_i = \min(\relw_i, 0)$.

In the following, we investigate the relation between the behavior of the original capital 
process ${\cal K}^{\cal P}_n$ and the one-sided capital process ${\cal K}^{{\cal P}^+}_n$.
If $M_n$ or $\relw_n$ changes the sign only finite number of times, then
the behavior of ${\cal K}_n^{{\cal P}^+}$ is fairly trivial. 
If $\relw_n$ is eventually all non-negative, then 
${\cal K}^{\cal P}_n$ and ${\cal K}_n^{{\cal P}^+}$ are asymptotically
equivalent.  On the other hand if $\relw_n$ is eventually negative,
then ${\cal K}_n^{{\cal P}^+}$ stays constant, whereas  the behavior
of ${\cal K}^{\cal P}_n$  is described in Theorem \ref{thm:rate}.
Therefore we will consider the case that $\relw_n$ changes the sign 
infinitely often.  At first, by writing 
\begin{align*}
\tilde{s}_n = (1 - \rho)\alpha - \rho \beta + s_n,
\end{align*}
we can express
\begin{align*}
\relw_n = \frac{\tilde{s}_{n-1}}{(\alpha + \beta + n - 1)\rho(1 - \rho)}.
\end{align*} 
% and recall that for the Bayesian strategy ${\cal P}$,
% \begin{align*}
% {\cal K}_n^{\cal P} = ({\rnp})^{-h_n}({(1-\rnp)})^{-t_n}
% \frac{(\alpha)_{h_n}(\beta)_{t_n}}{(\alpha + \beta)_n}.
% \end{align*}
Noting that at time $n$ when $\tilde{s}_n$ changes the sign,
\[ \frac{h_n}{n} = \rnp + O\Bigl(\frac{1}{n}\Bigr), \]
we start at sufficiently large time $n_0$ such that $h_{n_0}/n_0 
\simeq \rnp$, and proceed to divide the sequence $\{\tilde{s}_n\}$ 
into the following two types of blocks. 
For $n_0 \le k + 1 \le l$, consider a block $\{k + 1, \dots , l\}$.
We call it a {\it non-negative block} if
\[ \tilde{s}_k < 0, \tilde{s}_{k+1} \ge 0, \tilde{s}_{k+2} \ge 0, \dots , 
\tilde{s}_l \ge 0, \tilde{s}_{l+1} < 0. \]
%holds, we call $\{k + 1, \dots , l\}$ a non-negative block, 
Similarly we call it a {\it negative block} if 
\[ \tilde{s}_k \ge 0, \tilde{s}_{k+1} < 0, \tilde{s}_{k+2} < 0, \dots , 
\tilde{s}_l < 0, \tilde{s}_{l+1} \ge 0.  \] 
% holds, we call $\{k + 1, \dots , l\}$ a negative block. 
By definition, negative and non-negative blocks are alternating.

We first consider a particularly simple case of the fair-coin game
$\rho = 1/2$ and symmetric prior $\alpha=\beta$. 
In this case $\tilde s_n = s_n =0$ when $\tilde s_n$ changes the sign,  and the lengths of the blocks are always even numbers.
Then for each block 
$\{k+1,\dots,l\}$ we have $h_k = t_k$, $h_l = t_l$ 
and the capital ratio
${\cal K}_l^{\cal P}/{\cal K}_k^{\cal P}$
is expressed as
\begin{align*}
\frac{{\cal K}_l^{\cal P}}{{\cal K}_k^{\cal P}} = 
2^{2m}\frac{(\alpha + h_k)_m^2}{(2\alpha + 2h_k)_{2m}} 
= \prod_{j=0}^{m-1} \frac{2(\alpha + h_k + j)}{2(\alpha + h_k + j) + 1} 
< 1,
\end{align*}
where $l - k = 2m$ is an even number. 
As for the one-sided capital ratio, 
${\cal K}^{{\cal P}^+}_l/{\cal K}^{{\cal P}^+}_k = 
{\cal K}_l^{\cal P}/{\cal K}_k^{\cal P}$  
during a non-negative block, and 
${\cal K}^{{\cal P}^+}_l/{\cal K}^{{\cal P}^+}_k = 1$ 
during a negative block. 
Therefore ${\cal K}_n^{\cal P} < {\cal K}^{{\cal P}^+}_n$ holds for
all $n$ in the non-negative block.

Now we consider a general biased-coin game.  When $\tilde s_n$ changes
the sign, we have to consider overshoot of order $O(1/n)$. 
Therefore we need to carefully bound the capital ratio 
${\cal K}_l^{\cal P}/{\cal K}_k^{\cal P}$ of  
${\cal P}$ and the capital ratio 
${\cal K}^{{\cal P}^+}_l/{\cal K}^{{\cal P}^+}_k$ of 
${\cal P}^+$ from above and below for each negative or non-negative block. 
This is conducted based on the log capital formula given by 
(\ref{eq:log capital}), that is 
\begin{align*}
\log {\cal K}_n^{\cal P}
= n' D\Bigl(\frac{h_n'}{n'} \Big\| \rnp\Bigr)  
-\frac{1}{2}\log\frac{ h_n' t_n'}{n'} +
c_0(\alpha,\beta) + O\biggl(\frac{1}{\min (h_n', t_n')}\biggr).
\end{align*}
In the above, when $n = k$ or $l$, 
\begin{align*}
\frac{h_n'}{n'} = \rnp + O\Bigl(\frac{1}{n} \Bigr),\qquad  
\frac{t_n'}{n'} = 1-\rnp + O\Bigl(\frac{1}{n} \Bigr), 
\end{align*}
and when $\delta$ is small, as was noted in Section \ref{sec:process},
\begin{align*}
D(\rnp + \delta \| \rnp) 
% = (\rnp + \delta) 
% \log \Bigl(1 + \frac{\delta}{\rnp}\Bigr) + 
% ((1-\rnp) - \delta) \log \Bigl(1 - \frac{\delta}{(1-\rnp)}\Bigr) 
= \frac{\delta^2}{2\rnp(1-\rnp)} + O(\delta^3).
\end{align*}
Therefore with  $\delta = \bar x_n= O(1/n)$ we have
\[
n' D\Bigl(\frac{h_n'}{n'} \Big\| \rnp\Bigr)  = O\Bigl(\frac{1}{n}\Bigr).
\]
Hence we get for $n = k$ or $l$
\begin{align*}
\log {\cal K}_n^{\cal P} = - \frac{1}{2}\log \rnp(1-\rnp) - \frac{1}{2}\log n + 
c_0(\alpha,\beta) + O\Bigl(\frac{1}{n} \Bigr),
\end{align*}
so that it follows
\begin{align*}
\log \frac{{\cal K}_l^{\cal P}}{{\cal K}_k^{\cal P}} = 
- \frac{1}{2}\log \frac{l}{k} + O\Bigl(\frac{1}{k} \Bigr),
\end{align*}
which implies that there exists a constant $C>0$ depending only on 
$\alpha, \beta,$ and $\rnp$ such that
\begin{align*}
- \frac{1}{2}\log \frac{l}{k} - \frac{C}{k} < 
\log \frac{{\cal K}_l^{\cal P}}{{\cal K}_k^{\cal P}} < 
- \frac{1}{2}\log \frac{l}{k} + \frac{C}{k}.
\end{align*}
We apply this relation for successive non-negative and negative blocks 
by noting the approximation formula
\begin{align*}
\frac{1}{m} + \frac{1}{m + 1} + \cdots + \frac{1}{n} \le 
\int_m^n \frac{dx}{x} + \frac{1}{m} = 
\log \frac{n}{m} + \frac{1}{m}.
\end{align*}
Then we obtain at an end point $n_l$ of any block,
\begin{align*}
\Bigl(- \frac{1}{2} - C \Bigr) \log \frac{n_l}{n_0} - \frac{C}{n_0} <
\log \frac{{\cal K}^{{\cal P}}_{n_l}}{{\cal K}^{{\cal P}}_{n_0}} < 
\Bigl(- \frac{1}{2} + C \Bigr) \log \frac{n_l}{n_0} + \frac{C}{n_0}.
\end{align*}

When we reach an end point $n_l$ after passing sufficiently many blocks,  
the one-sided capital ratios 
${\cal K}^{{\cal P}^+}_l/{\cal K}^{{\cal P}^+}_k$ 
behave in the same way as 
${\cal K}^{{\cal P}}_l/{\cal K}^{{\cal P}}_k$ 
during non-negative blocks and stay one during negative blocks, so that we can also bound 
$\log {\cal K}^{{\cal P}^+}_{n_l}/{\cal K}^{{\cal P}^+}_{n_0}$ 
in the following manner.
\begin{align*}
\Bigl(- \frac{1}{2} - C \Bigr) \log \frac{n_l}{n_0} - \frac{C}{n_0}  
< \log \frac{{\cal K}^{{\cal P}^+}_{n_l}}{{\cal K}^{{\cal P}^+}_{n_0}} < 
C\log \frac{n_l}{n_0} + \frac{C}{n_0}.
\end{align*}

{}From the above two relations, it follows that at a sufficiently large 
end point $n_l$, the one-sided log capital $\log {\cal K}_{n_l}^{{\cal P}^+}$ 
differs from the original log capital $\log {\cal K}_{n_l}^{{\cal P}}$ at most 
$O(\log n_l)$. 
This fact implies that the two log capitals $\log {\cal K}_n^{{\cal P}^+}$ and 
$\log {\cal K}_n^{\cal P}$ behave similarly except for $O(\log n)$. Thus from 
Theorem \ref{thm:rate}, we obtain the following result.

\begin{theorem}
\label{thm:one-sided}
The one-sided Bayesian strategy ${\cal P}^+$ given by 
(\ref{eq:one-sided strategy}) weakly forces the one-sided version of SLLN 
with the convergence rate $O(\sqrt{\log n/n})$  
and with the convergence factor $\sqrt{\rho(1-\rnp)}$, that is 
\begin{equation}
\limsup_n \frac{\sqrt{n}\bar{x}_n}{\sqrt{\log n}} > 
\sqrt{\rho(1-\rnp)}\  
\Rightarrow \ 
\limsup_n {\cal K}_n^{{\cal P}^+} = \infty.
\end{equation}
\end{theorem}

\section{Some numerical examples}
\label{sec:numerical}

In this section we illustrate capital processes of our strategies for
two data sets.  The first data set is based on the Nikkei 225 stock
average index for 500 days starting January 2000.  We set $x_n > 0$ if
the opening price of the $n+1$-st day was higher than the opening
price of the $n$-th day.  The second data set is based on the first
500 digits in the fractional part of $\pi-3=0.141592653589793...$.
We set $x_n > 0$ if the $n$-th digit is in $\{5,\ldots,9\}$.

For the Nikkei data we have $h_{500}=221$ heads and $t_{500}=279$
tails for the 500 days.  The values of $\log {\cal K}_{500}$ for the
hypergeometric model (HG) and the beta-binomial strategies for 
$\rho=1/2,\ 2/3,\ 2/5$, $\alpha=\beta=1,\ 100,\ 500$, together with the positive-part (PP)
and negative-part (NP) strategies, are tabulated in Table 1. %\ref{tab:1}.
We see the exponential growth of the capital process for $\rho = 2/3$.
For this example it seems to to be advantageous to take  $\alpha=\beta=100$.

\begin{table}[h]
\label{tab:1}
\caption{Log capital process at $n=500$ for Nikkei 225}
\begin{center}
\begin{tabular}{llll}\hline
$\alpha, \ \beta$ & $\rho=1/2$ & $\rho=2/3$ & $\rho=2/5$\\ \hline
HG           &  6.698416 &  56.24544  & 5.145427\\
1.0, 1.0     & 0.4818099 & 50.02884   & -1.071180\\ 
\quad  (PP)    & -1.712525 &0.0         & -1.071180\\
\quad  (NP)    & 2.194335  & 50.02884   & 0.0 \\
100, 100 & 1.781788  & 51.32881   & 0.2287981\\ 
\quad  (PP)    & -0.1556388& 0.0        & 0.2287981\\
\quad  (NP)    & 1.937426  & 51.32881   & 0.0 \\
500,  500 & 0.9195455 & 50.46657    & -0.633444 \\ 
\quad  (PP)    & -0.03559586 &0.0       & -0.633444\\
\quad  (NP)    & 0.9551413&50.46657     & 0.0
\end{tabular}
\end{center}
\end{table}

\begin{figure}[h]
\begin{center}
\includegraphics[width=14cm]{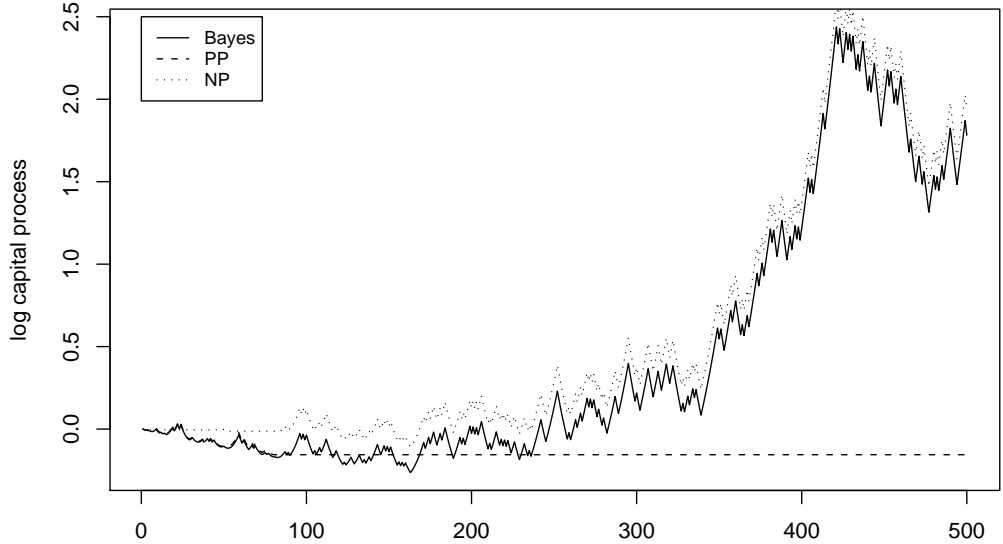}
\vspace*{-10mm}
\caption{Nikkei 225 data
($\rho=1/2,\ \alpha=\beta=100$)
}
\label{fig:nikkei}
\end{center}
\end{figure}

We plot the entire log capital processes of the beta-binomial strategy with its positive part and negative part for the case of $\rho=1/2,\ 
\alpha=\beta=100$ in Figure  \ref{fig:nikkei}. The log capital process of the beta-binomial strategy is plotted with a solid line, that of the positive-part strategy is plotted with a dashed line, and that of the negative-part strategy is plotted with a dotted line. We see that for this example the beta-binomial strategy is close to the negative-part strategy. The log capital of the positive-part strategy stays constant after about $n=90$.

For the digits of $\pi$, we have $h_{500}=239$ heads and $t_{500}=261$
tails.  The number of heads and tails are more balanced for this data
set than the Nikkei 225 case above.
Table 2 gives the same information as in Table 1 for this data
set.  We see the same tendency in Table 2 as in Table 1, although 
$\alpha=\beta=500$ performs better than $\alpha=\beta=100$.

\begin{table}[h]
\label{tab:2}
\caption{Log capital process at $n=500$ for the digits of $\pi$}
\begin{center}

\begin{tabular}{llll}\hline
$\alpha, \ \beta$ & $\rho=1/2$ & $\rho=2/3$ & $\rho=2/5$\\ \hline
HG           &  3.816784  & 40.88716  & 9.562166\\
1.0, 1.0     & -2.399822  & 34.67056  &  3.345560\\ 
\quad  (PP)  & -0.9942046 & 0.0      & 3.957064 \\
\quad  (NP)  & -1.4056175&  34.67056 & -0.6115032 \\
100,  100    & -0.2810085 & 36.78937 & 5.464374\\ 
\quad  (PP)  & -0.1820164 & 0.0      & 5.464374\\
\quad  (NP)  & -0.09899215 & 36.78937& 0.0 \\
500,  500 & -0.04136915 & 37.02901& 5.704013\\ 
\quad  (PP)  &  -0.04499401 & 0.0 &  5.704013\\  
\quad  (NP)  & 0.003624851 & 37.02901& 0.0
\end{tabular}
\end{center}
\end{table}

\begin{figure}[h]
\begin{center}
\includegraphics[width=14cm]{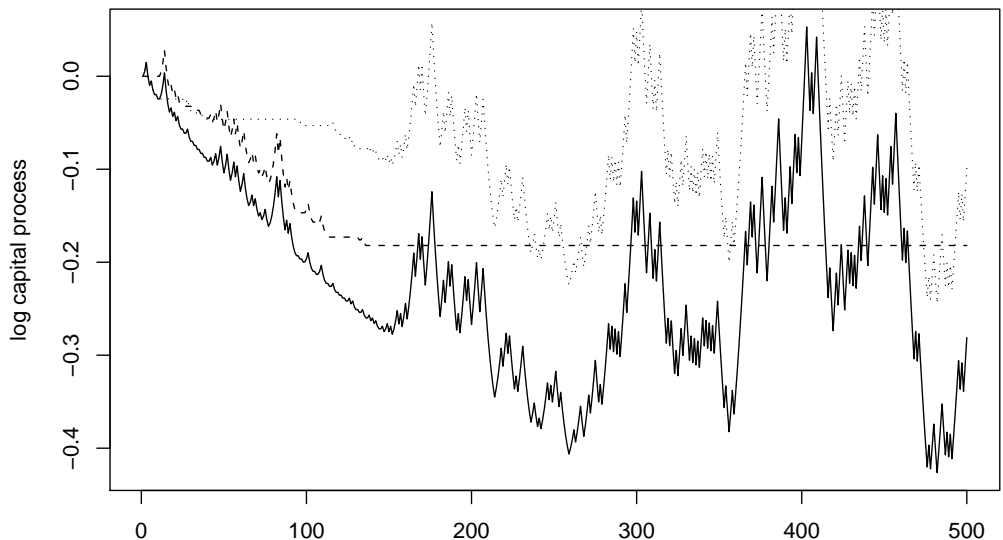}
\vspace*{-10mm}
\caption{500 digits of $\pi$\ ($\rho=1/2,\ \alpha=\beta=100$)
}
\label{fig:pi}
\end{center}
\end{figure}

As in Figure 1, the log capital processes of the beta-binomial strategy with its positive part and negative part for the case of $\rho=1/2,\ 
\alpha=\beta=100$ are plotted in Figure \ref{fig:pi}.  
The log capital of the positive-part strategy stays constant after about $n=140$.

\section{Concluding remarks}
\label{sec:remarks}

In this paper we have shown that for general biased-coin games, some simple Bayesian models provide explicit strategies of Skeptic which weakly force the strong law of large numbers with the convergence rate of
$O(\sqrt{\log n/n})$ and whose capital processes can be analyzed in
detail, leading naturally to the Kullback divergence.

We treated coin tossing and beta-binomial model for simplicity and
for the sake of explicit computations.  However we expect that many
of the results of this paper can be generalized to multidimensional
cases and more general prior distributions. Theorem \ref{thm:path}
should hold for games with unique risk neutral probability.

We have only considered exchangeable priors $Q$. Strategies based on
exchangeable priors can not exploit some block patterns or Markov
dependencies of Reality.  For example in the fair-coin game, the
beta-binomial prior can not exploit the following pattern of heads:
$1, 0, 1, 0,1, 0,\ldots$.  In order to exploit a variety of patterns
in Reality's moves, we need to use a prior $Q$ which contains
hyperparameters corresponding to these patterns.  One might consider
mixtures of priors, corresponding to various block patterns, Markov
dependence of various orders, etc.  
% Although exchangeable priors are
% probably too simple, it is not clear how to construct 
% practical priors which account for various higher order patterns.
We can take priors covering various higher order patterns and can  analyze optimalities with respect to such priors. Another interesting direction would be to extend our results to other limit theorems, such as the law of the iterated logarithm. These subjects will be treated in our subsequent works.

\end{document}